\documentclass{amsart}

\newtheorem{theorem}[equation]{Theorem}
\newtheorem{lemma}[equation]{Lemma}
\newtheorem{corollary}[equation]{Corollary}
\newtheorem{proposition}[equation]{Proposition}

\numberwithin{equation}{section}

\newcommand{\updot}{{\textstyle\cdot}}

\usepackage{amscd}

\begin{document}

\title{A cohomological property of Lagrange multipliers}
\author{Alan Adolphson}
\address{Department of Mathematics\\
Oklahoma State University\\
Stillwater, Oklahoma 74078}
\email{adolphs@math.okstate.edu}
\thanks{The first author was supported in part by NSA Grant
  \#MDA904-97-1-0068} 
\author{Steven Sperber}
\address{School of Mathematics\\
University of Minnesota\\
Minneapolis, Minnesota 55455}
\email{sperber@math.umn.edu}
\date{}
\begin{abstract}
The method of Lagrange multipliers relates the critical points of a
given function $f$ to the critical points of an auxiliary function
$F$.  We establish a cohomological relationship between $f$ and $F$
and use it, in conjunction with the Eagon-Northcott complex, to
compute the sum of the Milnor numbers of the critical points in
certain situations.
\end{abstract}
\maketitle

\section{Introduction}

Let $U$ be an open subset of ${\bf R}^n$ and let
$f,f_1,\ldots,f_r:U\rightarrow{\bf R}$ be continuously differentiable
functions on $U$.  Let $Y\subseteq U$ be the subset defined by
$f_1=\cdots=f_r=0$ and suppose that the matrix $(\partial f_i/\partial
x_j)_{\substack{i=1,\ldots,r\\ j=1,\ldots,n}}$ has rank $r$ at every
point of $Y$.  The usual theorem of Lagrange multipliers says that
${\bf a}=(a_1,\ldots,a_n)\in Y$ is a critical point of $f|_Y$ if and
only if there exists ${\bf b}=(b_1,\ldots,b_r)\in{\bf R}^r$ such that
$({\bf a};{\bf b})\in U\times{\bf R}^r$ is a critical point of the
auxiliary function $F=f+\sum_{i=1}^r y_if_i:U\times{\bf
  R}^r\rightarrow{\bf R}$.  The point ${\bf b}$ is unique when it
exists.

We establish a closer relation between $f$ and $F$ for algebraic
varieties over an arbitrary field $K$.  Let $X={\rm Spec}(A)$ be a
smooth affine $K$-scheme of finite type, purely of dimension $n$, let
$f,f_1,\ldots,f_r\in A$ and put $I=(f_1,\ldots,f_r)\subseteq A$.  Put
$B=A/I$ and let $Y={\rm Spec}(B)$, a closed subscheme of $X$.  We
assume that $Y$ is a smooth $K$-scheme, purely of codimension $r$ in
$X$.  We write $\bar{f}$ for the image of
$f\in A$ under the natural map $A\rightarrow B$.  Let $y_1,\ldots,y_r$
be indeterminates and consider $X\times_K{\bf A}^r={\rm
  Spec}(A[y_1,\ldots,y_r])$.  We shall write $A[y]$ for
$A[y_1,\ldots,y_r]$.  Put $F=f+\sum_{i=1}^r y_if_i\in
A[y]$.  Let $\Omega^k_{B/K}$ (resp.\ $\Omega^k_{A[y]/K}$) be the
module of differential $k$-forms of $B$ (resp.\ $A[y]$) over $K$.  Let
$d_{B/K}\bar{f}\in\Omega^1_{B/K}$ and $d_{A[y]/K}F\in\Omega^1_{A[y]/K}$ be the
exterior derivatives of $\bar{f}$ and $F$, respectively.  We consider the
complexes $(\Omega^{\updot}_{B/K},\phi_{\bar{f}})$ and
$(\Omega^{\updot}_{A[y]/K},\phi_F)$, where
$\phi_{\bar{f}}:\Omega^k_{B/K}\rightarrow\Omega^{k+1}_{B/K}$
is the map defined by
\[ \phi_{\bar{f}}(\omega)=d_{B/K}\bar{f}\wedge\omega \]
and $\phi_F:\Omega^k_{A[y]/K}\rightarrow\Omega^{k+1}_{A[y]/K}$
is the map defined by
\[ \phi_F(\omega)=d_{A[y]/K}F\wedge\omega. \]
The cohomology of these complexes is supported on the sets of critical
points of $\bar{f}$ and $F$, respectively.  The purpose of this note
is to prove the following result.  

\begin{theorem}
With the above notation and hypotheses, there are $A$-module
isomorphisms for all $i$
\[ H^i(\Omega^{\updot}_{B/K},\phi_{\bar{f}})\simeq
H^{i+2r}(\Omega^{\updot}_{A[y]/K},\phi_F). \]
\end{theorem}

The cohomology group $H^{n-r}(\Omega^{\updot}_{B/K},\phi_{\bar{f}})$
plays an important role.  For example, if $\bar{f}$ has only isolated
critical points on $Y$, then
$H^{n-r}(\Omega^{\updot}_{B/K},\phi_{\bar{f}})$ is a
finite-dimensional $K$-vector space whose dimension equals the sum of
the Milnor numbers of the critical points of $\bar{f}$ on $Y$.  In this
case, $H^i(\Omega^{\updot}_{B/K},\phi_{\bar{f}})=0$ for all $i\neq
n-r$.  To see this, since the assertion is local, we may assume that
$\Omega^1_{B/K}$ is a free $B$-module of rank $n-r$.  We may then
choose derivations $D_1,\ldots,D_{n-r}\in{\rm Der}_K(B)$
that form a basis for ${\rm Der}_K(B)$ as $B$-module and identify
$(\Omega^{\updot}_{B/K},\phi_{\bar{f}})$ with the cohomological Koszul
complex on $B$ defined by $D_1\bar{f},\ldots,D_{n-r}\bar{f}$.  In particular, 
\[ H^{n-r}(\Omega^{\updot}_{B/K},\phi_{\bar{f}})\simeq
B/(D_1\bar{f},\ldots,D_{n-r}\bar{f}). \]
Since this is assumed to be finite-dimensional, the ideal
$(D_1\bar{f},\ldots,D_{n-r}\bar{f})$ of $B$ has height $n-r$,
therefore has depth $n-r$ as well (${\rm Spec}(B)$ smooth implies in
particular that $B$ is Cohen-Macaulay).  The depth sensitivity of the
Koszul complex \cite[Theorem 16.8 and Corollary]{MA} then implies that
all its cohomology in degree $<n-r$ vanishes (and that
$D_1\bar{f},\ldots,D_{n-r}\bar{f}$ form a regular sequence in $B$).

Let $\Omega^k_{A/K}$ be the module of differential $k$-forms of $A$
over $K$ and let $d_{A/K}:\Omega^k_{A/K}\rightarrow\Omega^{k+1}_{A/K}$
be exterior differentiation.  As a corollary of the proof of
Theorem~1.1, we shall obtain the following. 
\begin{theorem}
Under the hypothesis of Theorem $1.1$, there is an isomorphism of
$B$-modules 
\[ H^{n-r}(\Omega^{\updot}_{B/K},\phi_{\bar{f}})\simeq 
\biggl(\Omega^n_{A/K}\bigg/(d_{A/K}f\wedge d_{A/K}f_1\wedge\cdots\wedge
d_{A/K}f_r\wedge\Omega^{n-r-1}_{A/K})\biggr)\bigotimes_A B. \]
\end{theorem}

We write out this isomorphism in a special case.  Let $X={\bf A}^n$,
so that $A$ is the polynomial ring $K[x_1,\ldots,x_n]$ and
$f,f_1,\ldots,f_r\in K[x_1,\ldots,x_n]$.  Let $J$ be the ideal of
$K[x_1,\ldots,x_n]$ generated by the $(r+1)\times(r+1)$-minors of the
matrix 
\begin{equation} \begin{bmatrix}
\partial f_1/\partial x_1 & \dots & \partial f_1/\partial x_n \\
\hdotsfor{3} \\
\partial f_r/\partial x_1 & \dots & \partial f_r/\partial x_n \\
\partial f/\partial x_1 & \dots & \partial f/\partial x_n 
\end{bmatrix}.
\end{equation}
\begin{corollary}
If $f_1,\ldots,f_r\in K[x_1,\ldots,x_n]$ define a smooth complete
intersection in ${\bf A}^n$, then
\[ H^{n-r}(\Omega^{\updot}_{B/K},\phi_{\bar{f}})\simeq
K[x_1,\ldots,x_n]/(I+J). \]
\end{corollary}

Again assume $f,f_1,\ldots,f_r\in K[x_1,\ldots,x_n]$.  Let
$d_i=\deg f_i$ for $i=1,\ldots,r$, $d_{r+1}=\deg f$, and denote by
$f_i^{(d_i)}$  (resp.\ $f^{(d_{r+1})}$) the homogeneous part of $f_i$
(resp.\ of $f$) of highest degree.  For $i=1,\ldots,r$, Let
$\tilde{f}_i$ be the homogenization of $f_i$, i.~e., 
\[ \tilde{f}_i=x_0^{d_i}f_i(x_1/x_0,\ldots,x_n/x_0)\in
K[x_0,\ldots,x_n]. \]
Using Corollary 1.4 and the Eagon-Northcott complex\cite{EN}, we shall
prove the following result, which was suggested by a theorem of
Katz\cite[Th\'{e}or\`{e}me 5.4.1]{KA}.
\begin{proposition}
Suppose that $\tilde{f}_1=\cdots=\tilde{f}_r=0$ defines a smooth
complete intersection in ${\bf P}^n$ that intersects the hyperplane
$x_0=0$ transversally and that
$f_1^{(d_1)}=\cdots=f_r^{(d_r)}=f^{(d_{r+1})}=0$ defines a smooth complete
intersection in ${\bf P}^{n-1}$.  If ${\rm char}(K)>0$, we assume also
that $(d_{r+1},{\rm char}(K))=1$.  Then $\bar{f}$ has only isolated critical
points on the variety $Y\subseteq{\bf A}^n$ defined by
$f_1=\cdots=f_r=0$ and $\dim_K
H^{n-r}(\Omega^{\updot}_{B/K},\phi_{\bar{f}})$ equals the coefficient
of $t^{n-r}$ in the power series expansion at $t=0$ of the rational
function
\[ \frac{d_1\cdots d_r(1-t)^n}{\prod_{i=1}^{r+1}(1-d_it)}. \]
\end{proposition}

\section{An intermediate complex}

We reduce Theorem 1.1 to a related statement.  Keeping our hypotheses
on $X$ and~$Y$, we shall express the complex
$(\Omega^{\updot}_{A[y]/K},\phi_F)$ as the
total complex associated to a certain double complex and show
that the vertical cohomology of this double complex vanishes except in
degree $r$.  The intermediate complex referred to in the title of this
section will be the horizontal complex associated to this single
nonvanishing vertical cohomology group.

We write $K[y]$ for $K[y_1,\ldots,y_r]$.  Let
\[ d_1:\Omega^k_{A[y]/K[y]}\rightarrow\Omega^{k+1}_{A[y]/K[y]} \]
and
\[ d_2:\Omega^k_{A[y]/A}\rightarrow\Omega^{k+1}_{A[y]/A} \]
be the exterior derivatives.  For $p,q\geq 0$ put
\[ C^{p,q}=\Omega^p_{A[y]/K[y]}\bigotimes_{A[y]}\Omega^q_{A[y]/A}. \]
We have
\begin{equation}
\Omega^k_{A[y]/K}\simeq \bigoplus_{p+q=k}C^{p,q}. 
\end{equation}
Define $\delta_1:C^{p,q}\rightarrow C^{p+1,q}$ and
$\delta_2:C^{p,q}\rightarrow C^{p,q+1}$ by
\begin{align*}
\delta_1(\omega\otimes\omega')& = (d_1F\wedge\omega)\otimes\omega' \\
\delta_2(\omega\otimes\omega')& =
(-1)^p(\omega\otimes(d_2F\wedge\omega')).
\end{align*}
It is straightforward to check from (2.1) and these definitions that
$(\Omega^{\updot}_{A[y]/K},\phi_F)$ is the total complex of the
double complex $\{C^{p,q}\}$. 

Let $\phi'_F:\Omega^k_{A[y]/A}\rightarrow\Omega^{k+1}_{A[y]/A}$ be
defined by
\[ \phi'_F(\omega')=d_2F\wedge\omega'. \]
In coordinate form, we have 
\[ d_2F=\sum_{i=1}^r \frac{\partial F}{\partial
  y_i}\,dy_i=\sum_{i=1}^r f_i\,dy_i, \]
so the complex $(\Omega^{\updot}_{A[y]/A},\phi'_F)$ is isomorphic
to the (cohomological) Koszul complex on $A[y]$ defined by
$f_1,\ldots,f_r$.  This Koszul complex decomposes into a direct sum of
copies (indexed by the monomials in $y_1,\ldots,y_r$) of the Koszul
complex on $A$ defined by $f_1,\ldots,f_r$.  We denote this latter
Koszul complex by ${\rm Kos}(A;f_1,\ldots,f_r)$.  Our hypothesis that
$Y$ is a smooth complete intersection defined by $f_1,\ldots,f_r$
implies that $f_1,\ldots,f_r$ form a regular sequence in the local
ring of $X$ at any point of $Y$.  This gives the following result.

\begin{lemma} 
The cohomology of the complex ${\rm Kos}(A;f_1,\ldots,f_r)$ is given
by
\begin{align*}
  H^i({\rm Kos}(A;f_1,\ldots,f_r)) & =0 \quad\text{if $i\neq r$}, \\
  H^r({\rm Kos}(A;f_1,\ldots,f_r)) & = B.
\end{align*}
\end{lemma}

In view of the remarks preceding the lemma, we get the following.
\begin{corollary}
The cohomology of the complex
$(\Omega^{\updot}_{A[y]/A},\phi'_F)$ is given by
\begin{align}
H^i(\Omega^{\updot}_{A[y]/A},\phi'_F)& =0 \quad\text{if $i\neq r$}, \\
H^r(\Omega^{\updot}_{A[y]/A},\phi'_F)& =B[y_1,\ldots,y_r].
\end{align}
\end{corollary}

Equation (2.4) says that the vertical cohomology of the double complex
$\{C^{p,q}\}$ vanishes except in degree $r$.  Equation (2.5) and
standard results in homological algebra relating the vertical
cohomology of a double complex to the cohomology of its associated
total complex then imply that
\begin{equation}
H^{i+r}(\Omega^{\updot}_{A[y]/K},\phi_F)\simeq H^i\biggl(
\Omega^{\updot}_{A[y]/K[y]}\bigotimes_{A[y]}B[y],\bar{\delta}_1\biggr) 
\end{equation}
for all $i$, where $\bar{\delta}_1$ is the map induced by $\delta_1$.

We have isomorphisms
\begin{equation}
\Omega^k_{A[y]/K[y]}\bigotimes_{A[y]}B[y] \simeq
\Omega^k_{A/K}\bigotimes_A B[y]
\end{equation}
for all $k$.  Let $d_{A/K}:\Omega^k_{A/K}\rightarrow\Omega^{k+1}_{A/K}$ be
the exterior derivative.  By abuse of notation, we denote by
$d_{A/K}F$ the element 
\begin{equation}
d_{A/K}F=d_{A/K}f\otimes 1 +\sum_{i=1}^r d_{A/K}f_i\otimes y_i \in
\Omega^1_{A/K}\bigotimes_A B[y]. 
\end{equation}
There is a canonical map
\[ \biggl(\Omega^k_{A/K}\bigotimes_A B[y]\biggr) \bigotimes_{B[y]}
\biggl(\Omega^l_{A/K}\bigotimes_A B[y]\biggr) \rightarrow
\Omega^{k+l}_{A/K}\bigotimes_A B[y] \]
that sends $(\omega_1\otimes\alpha_1)\otimes
(\omega_2\otimes\alpha_2)$ to $(\omega_1\wedge\omega_2)\otimes
(\alpha_1\alpha_2)$.  We denote by 
\[ (\omega_1\otimes\alpha_1)\wedge (\omega_2\otimes\alpha_2) \]
the image of $(\omega_1\otimes\alpha_1)\otimes
(\omega_2\otimes\alpha_2)$ under this map.  Define
$\tilde{\phi}_F:\Omega^k_{A/K}\bigotimes_A B[y]\rightarrow 
\Omega^{k+1}_{A/K}\bigotimes_A B[y]$ by
\[ \tilde{\phi}_F(\omega\otimes\alpha)=d_{A/K}F\wedge(\omega\otimes\alpha). \]
It is straightforward to check that under the isomorphism (2.7), the 
map $\bar{\delta}_1$ on the left-hand side is identified with the map
$\tilde{\phi}_F$ on the right-hand side.  Thus (2.6) gives
isomorphisms for all $i$
\begin{equation}
H^{i+r}(\Omega^{\updot}_{A[y]/K},\phi_F)\simeq
H^i(\Omega^{\updot}_{A/K}\bigotimes_A B[y],\tilde{\phi}_F).
\end{equation}
The main effort in this paper will be devoted to proving the
following.

\begin{theorem}
There is an injective quasi-isomorphism of complexes of $A$-modules
\[ (\Omega^{\updot}_{B/K}[-r],\phi_{\bar{f}})\hookrightarrow
(\Omega^{\updot}_{A/K}\bigotimes_A B[y],\tilde{\phi}_F), \]
where $\Omega^{\updot}_{B/K}[-r]$ is the complex with
$\Omega^i_{B/K}[-r]=\Omega^{i-r}_{B/K}$.  In particular, there are
isomorphisms for all $i$
\[ H^{i-r}(\Omega^{\updot}_{B/K},\phi_{\bar{f}})\simeq
H^i(\Omega^{\updot}_{A/K}\bigotimes_A B[y],\tilde{\phi}_F). \]
\end{theorem}

Theorem 1.1 clearly follows from Theorem 2.10 and equation (2.9).  The
proof of Theorem 2.10 will be carried out in sections 3 and 4.  In
section 3, we define an isomorphism of
$(\Omega^{\updot}_{B/K}[-r],\phi_{\bar{f}})$ with a subcomplex
$L^{\updot}$ of $(\Omega^{\updot}_{A/K} \bigotimes_A
B[y],\tilde{\phi}_F)$.  In section 4, we show that the inclusion of
$L^{\updot}$ into $\Omega^{\updot}_{A/K} \bigotimes_A B[y]$ is a
quasi-isomorphism.

\section{Proof of Theorem 2.10}

We begin by using the degree in the $y$-variables to construct an
increasing filtration on the complex
$(\Omega^{\updot}_{A/K}\bigotimes_A B[y],\tilde{\phi}_F)$.  Let $F_dB[y]$
be the $B$-module of all polynomials of degree $\leq d$ in
$y_1,\ldots,y_r$.  We define the filtration $F.$ on
$\Omega^{\updot}_{A/K}\bigotimes_A B[y]$ by setting
\[ F_d(\Omega^k_{A/K}\bigotimes_A B[y])=\Omega^k_{A/K}\bigotimes_A
F_{d+k}B[y]. \]
By (2.8) we see that
\[ \tilde{\phi}_F\biggl(F_d(\Omega^k_{A/K}\bigotimes_A
B[y])\biggr)\subseteq F_d(\Omega^{k+1}_{A/K}\bigotimes_A B[y]), \]
hence we have a filtered complex.  Since $F_dB[y]=0$ for $d<0$, we
have
\[ F_d(\Omega^k_{A/K}\bigotimes_A B[y])=0\quad\text{for $d<-k$.} \]
Furthermore, $F_0B[y]=B$, so we make the identification
\begin{equation}
F_{-k}(\Omega^k_{A/K}\bigotimes_A B[y])=\Omega^k_{A/K}\bigotimes_A B.
\end{equation}

We define a map
\[ \Phi:\Omega^k_{A/K}\bigotimes_A B\rightarrow
\Omega^{k+r}_{A/K}\bigotimes_A B \]
by the formula
\[ \Phi(\xi)=(-1)^{kr}(d_{A/K}f_1\otimes 1)\wedge\cdots\wedge
(d_{A/K}f_r\otimes 1)\wedge\xi \]
for $\xi\in\Omega^k_{A/K}\bigotimes_A B$.
\begin{lemma} 
$\ker\Phi=\sum_{i=1}^r (d_{A/K}f_i\otimes
1)\wedge(\Omega^{k-1}_{A/K}\bigotimes_A B)$. 
\end{lemma}

{\it Proof}.  It suffices to check equality locally, so we may assume
that $\Omega^1_{A/K}$ is a free $A$-module, of rank $n$.  Then
$\Omega^1_{A/K}\bigotimes_A B$ is a free $B$-module of rank $n$ and
\[ \Omega^k_{A/K}\bigotimes_A B\simeq
\bigwedge^k(\Omega^1_{A/K}\bigotimes_A B). \] 
We are thus in the
situation of Saito\cite{SA}.  The smooth complete intersection
hypothesis implies that the ideal of $B$ generated by the coefficients
of $(d_{A/K}f_1\otimes 1)\wedge\cdots\wedge (d_{A/K}f_r\otimes 1)$
relative to the basis of $\Omega^r_{A/K}\bigotimes_A B$ obtained by
taking $r$-fold exterior products of a basis of
$\Omega^1_{A/K}\bigotimes_A B$ is the unit ideal.  The desired
conclusion then follows from part (i) of the theorem of \cite{SA}.

It follows from Lemma 3.2 that
\[ (\Omega^k_{A/K}\bigotimes_A B)/\ker \Phi \simeq
\biggl(\Omega^k_{A/K}\bigg/\sum_{i=1}^r
(d_{A/K}f_i\wedge\Omega^{k-1}_{A/K})\biggr) \bigotimes_A B. \]
By a standard result, this is just $\Omega^k_{B/K}$.  
Using the identification (3.1) (with $k$ replaced by $k+r$), we see
that $\Phi$~induces an imbedding
\[ \bar{\Phi}:\Omega^k_{B/K}\hookrightarrow
\Omega^{k+r}_{A/K}\bigotimes_A B[y] \]
with image in $F_{-k-r}(\Omega^{k+r}_{A/K}\bigotimes_A B[y])$.
Equation (2.8) implies that
\begin{equation}
\tilde{\phi}_F({\Phi}(\xi))=(d_{A/K}f\otimes 1)\wedge{\Phi}(\xi) 
\end{equation}
for $\xi\in\Omega^k_{A/K}\bigotimes_A B$, from which it is easily
seen that $\bar{\Phi}$ is a map of complexes.

Equation (3.3) gives additional information.  Define a subcomplex
$(L^{\updot},\tilde{\phi}_F)$ of $(\Omega^{\updot}_{A/K}\bigotimes_A
B[y],\tilde{\phi}_F)$ by setting
\[ L^k=\{\xi\in F_{-k}(\Omega^k_{A/K}\bigotimes_A B[y])
\mid \tilde{\phi}_F(\xi)\in F_{-k-1}(\Omega^{k+1}_{A/K}\bigotimes_A
B[y])\}. \]
\begin{proposition}
The map $\bar{\Phi}$ is an isomorphism of complexes from
$(\Omega^{\updot}_{B/K}[-r],\phi_{\bar{f}})$ onto
$(L^{\updot},\tilde{\phi}_F)$. 
\end{proposition}

{\it Proof}.  Let $\omega\in\Omega^{k-r}_{B/K}$.  Then clearly
$\bar{\Phi}(\omega)\in F_{-k}(\Omega^k_{A/K}\bigotimes_A B[y])$ and,
by (3.3), $\bar{\Phi}(\omega)\in L^k$.  Thus $\bar{\Phi}$ is an
injective homomorphism of complexes whose image is contained in
$L^{\updot}$.  It only remains to prove that
$L^k\subseteq\bar{\Phi}(\Omega^{k-r}_{B/K})$ for all $k$.

Let $\xi\in F_{-k}(\Omega^k_{A/K}\bigotimes_A B[y])$.  We may
write
\[ \xi=\sum_j \omega_j\otimes b_j \]
with $\omega_j\in\Omega^k_{A/K}$ and $b_j\in B$.  We have
\[ d_{A/K}F\wedge\xi=\sum_j (d_{A/K}f\wedge\omega_j)\otimes b_j+\sum_{i=1}^r
\biggl(\sum_j(d_{A/K}f_i\wedge\omega_j)\otimes b_jy_i\biggr). \]
Since $\Omega^{k+1}_{A/K}\bigotimes_A B[y]$ is locally a free
$B[y]$-module, $\tilde{\phi}_F(\xi)\in
F_{-k-1}(\Omega^{k+1}_{A/K}\bigotimes_A B[y])$ if and only if
\[ \sum_j (d_{A/K}f_i\wedge\omega_j)\otimes b_j=0\quad\text{for
  $i=1,\ldots,r$,} \]
i.~e., $\xi\in L^k$ if and only if $(d_{A/K}f_i\otimes 1)\wedge\xi=0$ for
$i=1,\ldots,r$.  Thus the proof will be completed by the following
result.
\begin{lemma} Suppose $\xi\in\Omega^k_{A/K}\bigotimes_A B$
  satisfies
\[ (d_{A/K}f_i\otimes 1)\wedge\xi=0\quad\text{for $i=1,\ldots,r$.} \]
Then $\xi\in{\rm im}\;\Phi$.
\end{lemma}

{\it Proof}.  It suffices to check the condition locally, i.~e., to
show that for any maximal ideal $\bar{\bf m}$ of $B$, if
\begin{equation}
(d_{A/K}f_i\otimes 1)_{\bar{\bf m}}\wedge\xi_{\bar{\bf m}}=0\quad\text{for
  $i=1,\ldots,r$,} 
\end{equation}
then $\xi_{\bar{\bf m}}\in{\rm im}\;\Phi_{\bar{\bf m}}$.  Let ${\bf
    m}$ be the maximal ideal of $A$ corresponding to the maximal ideal
  $\bar{\bf m}$ of $B$.  The smooth complete intersection hypotheses
  implies that $(d_{A/K}f_1)_{{\bf m}},\ldots,(d_{A/K}f_r)_{{\bf m}}$
  can be extended to a basis of $(\Omega^1_{A/K})_{{\bf m}}$ as
  $A_{{\bf m}}$-module.  This implies that $\{(d_{A/K}f_i\otimes
  1)_{\bar{\bf m}}\}_{i=1}^r$ can be extended to a basis of
  $(\Omega^1_{A/K}\otimes_A B)_{\bar{\bf m}}$ as $B_{\bar{\bf
      m}}$-module.  Since
\[ (\Omega^k_{A/K}\otimes_A B)_{\bar{\bf m}}\simeq
\bigwedge^k(\Omega^1_{A/K}\otimes_A B)_{\bar{\bf m}}, \]
our $k$-form $\xi_{\bar{\bf m}}$ can be written as a sum of $k$-fold
exterior products of these basis elements multiplied by elements
of $B_{\bar{\bf m}}$.  But then (3.6) implies that $(d_{A/K}f_i\otimes
1)_{\bar{\bf m}}$ must appear in each of these 
$k$-fold exterior products, i.~e., 
\[ \xi_{\bar{\bf m}}=(d_{A/K}f_1\otimes 1)_{\bar{\bf m}}\wedge\cdots\wedge
(d_{A/K}f_r\otimes 1)_{\bar{\bf m}}\wedge\eta_{\bar{\bf m}} \]
for some $\eta_{\bar{\bf m}}\in (\Omega^{k-r}_{A/K}\otimes_A
B)_{\bar{\bf m}}$.  Hence $\xi_{\bar{\bf m}}\in{\rm im}\;\Phi_{\bar{\bf m}}$.  

{\it Proof of Theorem $1.2$}.  Proposition 3.4 implies that
\[ H^{n-r}(\Omega^{\updot}_{B/K},\phi_{\bar{f}})\simeq
L^n/\tilde{\phi}_F(L^{n-1}). \]
But from the definition
\[ L^n= \Omega^n_{A/K}\bigotimes_A B \]
and from Proposition 3.4
\[ L^{n-1}=\Phi(\Omega^{n-r-1}_{A/K}\bigotimes_A B). \]
Theorem 1.2 now follows from equation (3.3) and the definition of $\Phi$.
\section{Completion of the proof}

To complete the proof, we show that the inclusion
$L^{\updot}\hookrightarrow \Omega^{\updot}_{A/K}\bigotimes_A B[y]$ is
a quasi-isomorphism.  For this, it suffices to show that the
corresponding map of associated graded complexes relative to the
filtration $F.$ is a quasi-isomorphism.  We write ${\rm gr}_d$ to
denote these associated graded complexes.

It is easily seen that
\[ {\rm gr}_d(\Omega^k_{A/K}\bigotimes_A B[y])=\Omega^k_{A/K}
\bigotimes_A B[y]^{(d+k)}, \]
where $B[y]^{(d)}$ denotes the $B$-module of homogeneous polynomials
of degree $d$ in $y_1,\ldots,y_r$.  Furthermore, the differential
${\rm gr}(\tilde{\phi}_F)$ of this associated graded complex is given
by 
\[ {\rm gr}(\tilde{\phi}_F)(\xi)= \sum_{i=1}^r
(d_{A/K}f_i\otimes y_i)\wedge\xi \]
for $\xi\in \Omega^k_{A/K}\bigotimes_A B[y]^{(d+k)}$.  
It is also easy to see that $F.$ induces the ``stupid'' filtration on
$L^{\updot}$, i.~e.,
\[ F_dL^k=\begin{cases} L^k& \text{if $d\geq -k$}, \\ 0& \text{if
      $d<-k$}, \end{cases} \]
hence ${\rm gr}_d(L^{\updot})$ is the complex with $L^{-d}$ in
degree $-d$ and zeros elsewhere if $-n\leq d\leq 0$ and is the zero
complex otherwise.  Thus the assertion that
\[ {\rm gr}_d(L^{\updot})\hookrightarrow{\rm
  gr}_d(\Omega^{\updot}_{A/K} \bigotimes_A B[y]) \]
is a quasi-isomorphism is equivalent to the assertion that
\[ 0\rightarrow L^{-d}\rightarrow \Omega^{-d}_{A/K}\bigotimes_A
B[y]^{(0)}\rightarrow \Omega^{1-d}_{A/K}\bigotimes_A
B[y]^{(1)}\rightarrow\cdots\rightarrow \Omega^n_{A/K}\bigotimes_A
B[y]^{(d+n)}\rightarrow 0 \]
is exact for $-n\leq d\leq 0$ and
\[ 0\rightarrow \Omega^0_{A/K}\bigotimes_A
B[y]^{(d)}\rightarrow\cdots\rightarrow \Omega^n_{A/K}\bigotimes_A
B[y]^{(d+n)}\rightarrow 0 \]
is exact for $d>0$.  The definition of $L^{\updot}$ shows that the
sequence
\[ 0\rightarrow L^{-d}\rightarrow \Omega^{-d}_{A/K}\bigotimes_A
B[y]^{(0)}\rightarrow \Omega^{1-d}_{A/K}\bigotimes_A
B[y]^{(1)} \]
is exact for $-n\leq d\leq 0$.  Thus we need to show that the sequence
\[ \Omega^{k-1}_{A/K}\bigotimes_A B[y]^{(d+k-1)}\rightarrow
\Omega^k_{A/K}\bigotimes_A B[y]^{(d+k)}\rightarrow 
\Omega^{k+1}_{A/K}\bigotimes_A B[y]^{(d+k+1)} \]
is exact whenever $d>-k$.

This can be summarized in the following result.
\begin{proposition}
$H^k(\Omega^{\updot}_{A/K}\bigotimes_A B[y]^{(\cdot\,+d)},
{\rm gr}(\tilde{\phi}_F))=0$ for $d>-k$.
\end{proposition}

{\it Proof}.  The complex in question is a complex of $A$-modules
supported on ${\rm Spec}(B)$.  Hence to prove the desired vanishing,
we may first localize at a maximal ideal of $A$ containing
$f_1,\ldots,f_r$.  Thus we may assume that $A$ is a smooth local
$K$-algebra of dimension $n$ whose maximal ideal ${\bf m}$ contains
$f_1,\ldots,f_r$ and that $B=A/(f_1,\ldots,f_r)$ is a smooth local
$K$-algebra of dimension $n-r$.  This implies that there exist
$f_{r+1},\ldots,f_n\in{\bf m}$ such that
$d_{A/K}f_1,\ldots,d_{A/K}f_n$ form a basis for $\Omega^1_{A/K}$.
To simplify notation, we write
\[ \Omega_{i_1\cdots i_k}= (d_{A/K}f_{i_1}\otimes 1)\wedge\cdots\wedge
(d_{A/K}f_{i_k}\otimes 1)\in \Omega^k_{A/K}\bigotimes_A B[y]. \]
Then $\Omega^k_{A/K}\bigotimes_A B[y]$ is a free $B[y]$-module with
basis
\begin{equation}
\{\Omega_{i_1\cdots i_k} \mid 1\leq i_1<\cdots<i_k\leq n\}
\end{equation}
and is a free $B$-module with basis
\begin{equation}
\{y_1^{a_1}\cdots y_r^{a_r}\Omega_{i_1\cdots i_k} \mid
a_1,\ldots,a_r\in{\bf N}, \quad 1\leq i_1<\cdots<i_k\leq n\}.
\end{equation}
The differential ${\rm gr}(\tilde{\phi}_F)$ of the complex takes the
form
\[ {\rm gr}(\tilde{\phi}_F)(\xi)=\sum_{i=1}^r y_i\Omega_i\wedge\xi. \]

We proceed by induction on $r$.  Suppose $r=1$ and let
$\xi\in\Omega^k_{A/K}\bigotimes_A B[y_1]^{(d+k)}$.  The condition
$d>-k$ implies that $\xi$ is divisible by $y_1$, i.~e., $\xi$
can be written
\[ \xi=\sum_{1\leq i_1<\cdots<i_k\leq n}
y_1\xi(i_1,\ldots,i_k)\Omega_{i_1\cdots i_k} \]
with $\xi(i_1,\ldots,i_k)\in B[y_1]^{(d+k-1)}$.  The condition that
$\xi$ be a cocycle is that
\[ y_1\Omega_1\wedge\xi=0. \]
Since (4.2) is a basis, we see that this is the case if and only if
$\xi(i_1,\ldots,i_k)\neq 0$ implies $i_1=1$.  Put
\[ \eta=\sum_{2\leq i_2<\cdots<i_k\leq n} \xi(1,i_2,\ldots,i_k)
\Omega_{i_2\cdots i_k}. \]
Then $\eta\in\Omega^{k-1}_{A/K}\bigotimes_A B[y_1]^{(d+k-1)}$ and
\[ \xi=y_1\Omega_1\wedge\eta, \]
so $\xi$ is a coboundary.

Now let $r\geq 2$ and suppose the proposition is true for $r-1$.  Let
\[ \xi\in\Omega^k_{A/K}\bigotimes_A B[y]^{(d+k)} \]
be a cocycle, i.~e.,
\begin{equation}
\sum_{i=1}^r y_i\Omega_i\wedge\xi=0.
\end{equation}
Let $h$ be the highest power of $y_1$ appearing in any term of
$\xi$ (in the decomposition corresponding to the basis (4.3)) and let
$\xi^{(h)}$ be the sum of all terms of degree $h$ in $y_1$.
Suppose $h>0$.  Looking at the terms of degree $h+1$ in $y_1$ in
equation (4.4) gives
\[ \Omega_1\wedge\xi^{(h)}=0, \]
hence
\[ \xi^{(h)}=\sum_{2\leq i_2<\cdots<i_k\leq n}
y_1\eta(i_2,\ldots,i_k)\Omega_{1i_2\cdots i_k} \]
for some $\eta(i_2,\ldots,i_k)\in B[y]^{(d+k-1)}$.  Put
\[ \eta=\sum_{2\leq i_2<\cdots<i_k\leq n}
\eta(i_2,\ldots,i_k)\Omega_{i_2\cdots i_k}. \]
Then $\eta\in\Omega^{k-1}_{A/K}\bigotimes_A B[y]^{(d+k-1)}$ and the
highest power of $y_1$ appearing in
\[ \xi-\sum_{i=1}^r y_i\Omega_i\wedge\eta \]
is $\leq h-1$.

We may thus reduce to the case $h=0$, i.~e., $y_1$ does not appear in
$\xi$.  Equation~(4.4) then implies the two equalities
\begin{align}
\Omega_1\wedge\xi& =0 \\
\sum_{i=2}^r y_i\Omega_i\wedge\xi& =0.
\end{align}
From (4.5) we have
\[ \xi=\sum_{2\leq i_2<\cdots<i_k\leq n} \xi(i_2,\ldots,i_k)
\Omega_{1i_2\cdots i_k} \]
with $\xi(i_2,\ldots,i_k)\in B[y_2,\ldots,y_r]^{(d+k)}$.  Put
\[ \xi'=\sum_{2\leq i_2<\cdots<i_k\leq n}\xi(i_2,\ldots,i_k)
\Omega_{i_2\cdots i_k}\in\Omega^{k-1}_{A/K}\bigotimes_A
B[y_2,\ldots,y_r]^{(d+k)}. \]
By (4.6) we have
\[ \sum_{i=2}^r y_i\Omega_i\wedge\xi' =0, \]
i.~e., $\xi'$ is a $(k-1)$-cocycle in the complex
$(\Omega^{\updot}_{A/K}\bigotimes_A B[y_2,\ldots,y_r]^{(\cdot\,+d+1)},
\bar{\phi})$, where
\[ \bar{\phi}(\zeta)=\sum_{i=2}^r y_i\Omega_i\wedge\zeta. \]
But $d>-k$ implies $d+1>-(k-1)$, so the induction hypothesis for $r-1$
applies and we conclude that $\xi'$ is a coboundary.  This means there
exists
\[ \eta'\in\Omega^{k-2}_{A/K}\bigotimes_A B[y_2,\ldots,y_r]^{(d+k-1)}
\]
such that
\[ \sum_{i=2}^r y_i\Omega_i\wedge\eta'=\xi'. \]
If we put
$\eta=-\Omega_1\wedge\eta'\in\Omega^{k-1}_{A/K}\bigotimes_A
B[y]^{(d+k-1)}$, then
\begin{align*}
\sum_{i=1}^r y_i\Omega_i\wedge\eta& = \Omega_1\wedge\xi' \\
 &=\xi,
\end{align*}
hence $\xi$ is a coboundary.

\section{Proof of Proposition 1.5}

Let $f,f_1,\ldots,f_r\in K[x_1,\ldots,x_n]$ satisfy the hypotheses of
Proposition 1.5 and let $Y\subseteq{\bf A}^n$ be the variety
$f_1=\cdots=f_r=0$.  We begin by showing that $\bar{f}$ has only isolated
critical points on $Y$.  For notational convenience we write $K[x]$
for $K[x_1,\ldots,x_n]$.  Let $I\subseteq K[x]$ be the ideal generated
by $f_1,\ldots,f_r$ and let $J\subseteq K[x]$ be the ideal generated
by the $(r+1)\times(r+1)$-minors of the matrix (1.3).  
Put 
\[ Z=V(I+J)\subseteq{\bf A}^n. \]
The underlying point set of $Z$ is the set of critical points of $\bar{f}$
on $Y$.  We wish to show it is finite.  If not, then $\dim Z\geq 1$,
so $\dim\tilde{Z}\geq 1$ also, where $\tilde{Z}$ denotes the closure of
$Z$ in ${\bf P}^n$ under the natural compactification by adjoining the
hyperplane $x_0=0$ at infinity.  This would imply that the
intersection $\tilde{Z}\cap\{x_0=0\}$ is nonempty.  We prove that in
fact it is empty, therefore $Z$ must be finite.

Consider the matrix
\begin{equation} \begin{bmatrix}
\partial f_1^{(d_1)}/\partial x_1 & \dots & \partial
f_1^{(d_1)}/\partial x_n \\ 
\hdotsfor{3} \\
\partial f_r^{(d_r)}/\partial x_1 & \dots & \partial
f_r^{(d_r)}/\partial x_n \\ 
\partial f^{(d_{r+1})}/\partial x_1 & \dots & \partial
f^{(d_{r+1})}/\partial x_n  
\end{bmatrix}.
\end{equation}
Let $I'\subseteq K[x]$ be the ideal generated by
$f_1^{(d_1)},\ldots,f_r^{(d_r)}$ and let $J'$ be the ideal generated
by the $(r+1)\times(r+1)$-minors of the matrix (5.1).  For a
homogeneous ideal $Q\subseteq K[x]$, we denote by $V(Q)\subseteq{\bf
  A}^n$ the affine variety it defines and by $\tilde{V}(Q)\subseteq{\bf
  P}^{n-1}$ the projective variety it defines.  As point sets we have
\[ \tilde{Z}\cap\{x_0=0\}=\tilde{V}(I'+J'). \]
Suppose there exists a point $z\in\tilde{Z}\cap\{x_0=0\}$.  By hypothesis,
$f_1^{(d_1)}=\cdots=f_r^{(d_r)}=0$ defines a smooth complete
intersection in ${\bf P}^{n-1}$, so at any set of homogeneous
coordinates for the point $z$ the first $r$
rows of the matrix (5.1) are linearly independent.  But since all
$(r+1)\times(r+1)$-minors of (5.1) vanish at any set of homogeneous
coordinates for $z$, the last row must be a
linear combination of the first $r$ rows, say, 
\[ \biggl(\frac{\partial f^{(d_{r+1})}}{\partial
  x_1},\ldots,\frac{\partial f^{(d_{r+1})}}{\partial x_n}\biggr) =
  \sum_{i=1}^r c_i \biggl(\frac{\partial f_i^{(d_i)}}{\partial
  x_1},\ldots,\frac{\partial f_i^{(d_i)}}{\partial x_n}\biggr) \]
when evaluated at homogeneous coordinates for $z$, where the $c_i$ lie
  in the algebraic closure of $K$.  For $j=1,\ldots,n$, we then have 
\[ x_j\frac{\partial f^{(d_{r+1})}}{\partial x_j} = \sum_{i=1}^r
  c_ix_j\frac{\partial f_i^{(d_i)}}{\partial x_j} \]
when evaluated at homogeneous coordinates for $z$.  Summing these
equations over $j$ and using the Euler relation gives
\[ d_{r+1}f^{(d_{r+1})}=\sum_{i=1}^r c_id_if_i^{(d_i)} \]
when evaluated at these homogeneous coordinates.  Each
$f_i^{(d_i)}$ vanishes at $z$ and we are assuming $(d_{r+1},{\rm
  char}(K))=1$ if ${\rm char}(K)>0$, therefore $f^{(d_{r+1})}$
vanishes at $z$ also.  But this contradicts the
hypothesis that $f^{(d_{r+1})}=f_1^{(d_1)}=\cdots=f_r^{(d_r)}=0$
defines a smooth complete intersection in ${\bf P}^{n-1}$.  

To calculate $\dim_K H^{n-r}(\Omega^{\updot}_{B/K},\phi_{\bar{f}})$, we
begin by recalling the definition of the Eagon-Northcott
complex\cite{EN}.  In order to facilitate application of the results of
\cite{EN,EN2}, we describe these complexes homologically, rather than
cohomologically.  For $p=0,1,\ldots,n-r$, let $C^{(1)}_p$ be the free
$K[x]$-module with the following basis.  For $p=0$, $C^{(1)}_0=K[x]$
with basis $1$.  For $p\geq 1$, $C^{(1)}_p$ has basis
\[ \xi_{i_1}\cdots\xi_{i_{p+r}}\eta_1^{j_1}\cdots\eta_{r+1}^{j_{r+1}},
\]
where $1\leq i_1<\cdots<i_{p+r}\leq n$ and $j_1,\ldots,j_{r+1}$ are
nonnegative integers satisfying
\[ j_1+\cdots+j_{r+1}=p-1. \]
Thus $C^{(1)}_p$ has rank $\binom{n}{p+r}\binom{p+r-1}{r}$.  Let
$\delta_1:C^{(1)}_p\rightarrow C^{(1)}_{p-1}$ be the
$K[x]$-module homomorphism whose action on basis
elements is defined as follows.  For $p=1$,
\[ \delta_1(\xi_{i_1}\cdots\xi_{i_{r+1}})=\frac{\partial(f_1,
  \ldots,f_r,f)}{\partial(x_{i_1},\cdots,x_{i_{r+1}})}. \]
For $p>1$,
\begin{multline*}
\delta_1(\xi_{i_1}\cdots\xi_{i_{p+r}}\eta_1^{j_1}\cdots\eta_{r+1}^{j_{r+1}})=
\\
\sum_{\substack{l=1 \\ j_l>0}}^{r+1}\sum_{m=1}^{p+r}
(-1)^{m-1}\frac{\partial f_l}{\partial x_{i_m}} \xi_{i_1}\cdots
\hat{\xi}_{i_m}\cdots \xi_{i_{p+r}}\eta_1^{j_1}\cdots
\eta_l^{j_l-1}\cdots \eta_{r+1}^{j_{r+1}}. 
\end{multline*}
This is the Eagon-Northcott complex associated to the matrix (1.3).

We define a grading and filtration on $C_{\updot}^{(1)}$ as follows.
Let $K[x]$ have the usual grading and define 
\[ {\rm degree}(\xi_{i_1}\cdots\xi_{i_{p+r}}\eta_1^{j_1}\cdots
\eta_{r+1}^{j_{r+1}})=(j_1+1)d_1+\cdots+(j_{r+1}+1)d_{r+1}-(p+r). \]
This defines a grading on $C_p^{(1)}$ which in turn defines a
filtration on $C_p^{(1)}$ by letting level $k$ of the filtration be
the $K$-span of homogeneous elements of degree $\leq k$.  It is
straighforward to check that $\delta_1$ preserves this filtration.  We
let $(\bar{C}_{\updot}^{(1)},\bar{\delta}_1)$ be the associated graded
complex.  It is easily checked that
$(\bar{C}_{\updot}^{(1)},\bar{\delta}_1)$ is the Eagon-Northcott
complex associated to the matrix (5.1).  
\begin{proposition}
For $p>0$,
\[ H_p(\bar{C}_{\updot}^{(1)},\bar{\delta}_1)=0 \]
and the Hilbert-Poincar\'{e} series of the graded module
$H_0(\bar{C}_{\updot}^{(1)},\bar{\delta}_1)$ has the form
\[ \frac{G(t)(1-t)^{n-r}+H(t)(1-t)^{n-r+1}}{(1-t)^n},\]
where $G(t),H(t)$ are polynomials and $G(1)$ equals the coefficient of
$t^{n-r}$ in the power series expansion at $t=0$ of
\[ \frac{(1-t)^n}{\prod_{i=1}^{r+1} (1-d_it)}. \]
\end{proposition}

{\it Proof}.  We show that the ideal $J'$ has depth $n-r$.  The
vanishing of $H_p(\bar{C}_{\updot}^{(1)},\bar{\delta}_1)$ for $p>0$
then follows from \cite[Theorem 1]{EN} and the assertion about the
Hilbert-Poincar\'{e} series of
$H_0(\bar{C}_{\updot}^{(1)},\bar{\delta}_1)$ follows from Theorems 1,
2, and the Lemma in \cite{EN2}. (Note that, in the notation of
\cite{EN2}, we take $\mu_i=d_i$, $\nu_j=-1$, so that the entry in row
$i$, column $j$ of the matrix (5.1) is a homogeneous polynomial of
degree $\mu_i+\nu_j=d_i-1$.)  In fact, it is shown in \cite{EN} that
the depth of $J'$ is $\leq n-r$, so we only need to show that the
depth is $\geq n-r$.

Since the depth and the height of $J'$ are equal, it suffices to show
that the height of $J'$ is $\geq n-r$.  We proved at the beginning of
this section that $\tilde{V}(I'+J')$ is empty.  But this
says that $\tilde{V}(I')$ and $\tilde{V}(J')$ have empty intersection.  By
hypothesis, $\tilde{V}(I')$ is a smooth complete intersection, purely of
dimension $n-1-r$, hence all components of $\tilde{V}(J')$ must have
dimension $<r$.  This implies that all components of $V(J')$ in
${\bf A}^n$ have dimension $\leq r$, i.~e., have codimension $\geq
n-r$.  This is equivalent to the assertion that the height of $J'$ is
$\geq n-r$. 

Put $B'=K[x]/I'$.  We consider the related complex
$(\bar{C}_{\updot}^{(1)}\bigotimes_{K[x]}B',\bar{\delta}_1\otimes 1)$.  
\begin{proposition}
For $p>0$,
\[
H_p\biggl(\bar{C}_{\updot}^{(1)}\bigotimes_{K[x]}B',\bar{\delta}_1\otimes
1\biggr)=0. \] 
\end{proposition}

{\it Proof}.  The complex $\bar{C}_{\updot}^{(1)}\bigotimes_{K[x]}B'$
is the Eagon-Northcott complex of the image of the matrix (5.1) in
$B'$.  So by the same argument used in the proof of Proposition~5.2,
it suffices to show that the ideal $(J'+I')/I'$ in $B'$ has depth
$\geq n-r$.  Since $f_1^{(d_1)},\ldots,f_r^{(d_r)}$ is a regular
sequence, the ring $B'$ is Cohen-Macaulay.  Therefore the height and
the depth of $(J'+I')/I'$ are equal, so we are again reduced to
showing that the height of $(J'+I')/I'$ is $\geq n-r$.  Let ${\bf m}$
denote the ideal $(x_1,\ldots,x_n)$ of $K[x]$.  Since
$\tilde{V}(I'+J')=\emptyset$, the only prime ideal of $B'$ 
containing $(J'+I')/I'$ is ${\bf m}/I'$.  So it suffices to show that
${\bf m}/I'$ has height $\geq n-r$.  Let ${\bf p}$ be a minimal prime
ideal of $K[x]$ containing $I'$.  Since $\tilde{V}(I')$ is purely of
codimension $r$ in ${\bf P}^{n-1}$, ${\bf p}$ has height $r$.
Therefore every saturated chain of prime ideals from ${\bf p}$ to
${\bf m}$ has length $n-r$.  It follows that ${\bf m}/I'$ has height
$n-r$ in $B'$. 

Let $(C^{(2)}_{\updot},\delta_2)$ be the Koszul complex on
$K[x]$ defined by $f_1,\ldots,f_r$.  We regard $C^{(2)}_q$ as the free
$K[x]$-module with basis 
\[ \zeta_{k_1}\cdots\zeta_{k_q}, \]
where $1\leq k_1<\cdots<k_q\leq r$, and $\delta_2:C^{(2)}_q\rightarrow
C^{(2)}_{q-1}$ is defined by
\[ \delta_2(\zeta_{k_1}\cdots\zeta_{k_q})=\sum_{l=1}^q
(-1)^{l-1}f_{k_l}\,
\zeta_{k_1}\cdots\hat{\zeta}_{k_l}\cdots\zeta_{k_q}. \]
Each $C^{(2)}_q$ is graded by using the usual grading on
$K[x]$ and by defining the degree of
$\zeta_{k_1}\cdots\zeta_{k_q}$ to be $d_{k_1}+\cdots+d_{k_q}$.  The
complex $(C^{(2)}_{\updot},\delta_2)$ becomes a filtered complex by
defining level $k$ of the filtration to be the $K$-span of homogeneous
elements of degree $\leq k$.  Its associated graded complex
$(\bar{C}^{(2)}_{\updot},\bar{\delta}_2)$ is the Koszul complex on
$K[x]$ defined by $f_1^{(d_1)},\ldots,f_r^{(d_r)}$.  

Consider the double complex
$C^{(1)}_{\updot}\bigotimes_{K[x]} C^{(2)}_{\updot}$ and
let $T.$ be its associated total complex.  The filtrations on
$C^{(1)}_{\updot}$ and $C^{(2)}_{\updot}$ make $T_{\updot}$ a filtered
complex.  Its associated graded complex $\bar{T}.$ is the total
complex of the double complex $\bar{C}^{(1)}_{\updot}\bigotimes_{K[x]}
\bar{C}^{(2)}_{\updot}$.  One checks easily from the definitions that
\begin{align}
H_0(T.)& = K[x]/(I+J), \\
H_0(\bar{T}.)& = K[x]/(I'+J').
\end{align}
By (5.4) and Corollary 1.4,
\begin{equation}
H_0(T.)\simeq H^{n-r}(\Omega^{\updot}_{B/K},\phi_{\bar{f}}).
\end{equation}
We determine the relation between $H_0(T.)$ and $H_0(\bar{T}.)$.

Since $f_1^{(d_1)},\ldots,f_r^{(d_r)}$ is a regular sequence in $K[x]$, 
\[ H_p(\bar{C}^{(2)}_{\updot})=0\qquad\text{for $p>0$.} \]
Furthermore, $H_0(\bar{C}^{(2)}_{\updot})\simeq B'$.  This implies by
standard homological algebra that
\[ H_p(\bar{T}.)\simeq
H_p\biggl(\bar{C}^{(1)}_{\updot}\bigotimes_{K[x]} B'\biggr), \]
hence by Proposition 5.3
\begin{equation}
H_p(\bar{T}.)=0\qquad\text{for $p>0$.} 
\end{equation}
Standard homological algebra then implies that
\[ H_p(T.)=0\qquad\text{for $p>0$} \]
and that
\begin{equation}
{\rm gr}(H_0(T.))\simeq H_0(\bar{T}.) 
\end{equation}
as $K$-vector spaces, where the left-hand side denotes the associated
graded relative to the filtration induced by $T.$ on $H_0(T.)$.  In
particular, (5.6) implies that
\begin{equation}
\dim_K H^{n-r}(\Omega^{\updot}_{B/K},\phi_{\bar{f}}) = \dim_K
H_0(\bar{T}.)
\end{equation}

Proposition 5.2 and standard homological algebra imply that
\begin{equation}
H_p(\bar{T}.)\simeq H_p\biggl(H_0(\bar{C}^{(1)}_{\updot})\bigotimes_{K[x]}
\bar{C}^{(2)}_{\updot}\biggr). 
\end{equation}
Thus by (5.7), $H_0(\bar{C}^{(1)}_{\updot})\bigotimes_{K[x]}
\bar{C}^{(2)}_{\updot}$ is a resolution of
\begin{align*}
H_0\biggl(H_0(\bar{C}^{(1)}_{\updot})\bigotimes_{K[x]}
\bar{C}^{(2)}_{\updot}\biggr)& \simeq
H_0(\bar{C}^{(1)}_{\updot})\bigotimes_{K[x]}
H_0(\bar{C}^{(2)}_{\updot}) \\
& \simeq H_0(\bar{C}^{(1)}_{\updot})\bigotimes_{K[x]} B'.
\end{align*}
But $H_0(\bar{C}^{(1)}_{\updot})\bigotimes_{K[x]}
\bar{C}^{(2)}_{\updot}$is just the Koszul complex on
$H_0(\bar{C}^{(1)}_{\updot})$ defined by
$f_1^{(d_1)}$, \ldots, $f_r^{(d_r)}$.  It then follows from Proposition
5.2 that the Hilbert-Poincar\'{e} series of
$H_0(\bar{C}^{(1)}_{\updot})\bigotimes_{K[x]} B'$ is
\begin{equation}
\biggl(\prod_{i=1}^r
(1-t^{d_i})\biggr)\frac{G(t)(1-t)^{n-r}+H(t)(1-t)^{n-r+1}}{(1-t)^n}. 
\end{equation}
By (5.7) and (5.10), this is the Hilbert-Poincar\'{e} series of
$H_0(\bar{T}.)$, hence 
\[ \dim_K H_0(\bar{T}.)=d_1\cdots d_rG(1). \]
By Proposition 5.2, this proves that $\dim_K H_0(\bar{T}.)$ equals the
coefficient of $t^{n-r}$ in the power series expansion at $t=0$ of
\[ \frac{d_1\cdots d_r(1-t)^n}{\prod_{i=1}^{r+1}(1-d_it)}. \]
Proposition 1.5 then follows from (5.9).


\begin{thebibliography}{99}
\bibitem{AS} A. Adolphson and S. Sperber, Dwork cohomology, de Rham
  cohomology, and hypergeometric functions, Amer.\ J. Math.\ (to appear)
\bibitem{EN} J. A. Eagon and D. G. Northcott, Ideals defined by
  matrices and a certain complex associated with them, Proc.\ Roy.\
  Soc.\ A {\bf 269}(1962), 188--204
\bibitem{EN2} \bysame, A note on the Hilbert functions of certain
  ideals which are defined by matrices, Mathematika {\bf 9}(1962),
  118--126
\bibitem{KA} N. Katz, Sommes exponentielles, Ast\'{e}risque {\bf
    79}(1980), 1--209
\bibitem{MA} H. Matsumura, Commutative ring theory, Cambridge Studies
  in Advanced Mathematics {\bf 8}, Cambridge University Press, Cambridge,
  1986
\bibitem{SA} K. Saito, On a generalization of de Rham lemma, Ann.\
  Inst.\ Fourier, Grenoble {\bf 26}(1976), 165--170
\end{thebibliography}
\end{document}